\documentclass[12pt,a4paper]{article}
\usepackage{indentfirst}
\usepackage{amssymb}
\usepackage{amsmath}
\newtheorem{theo}{Theorem}[section]
\newtheorem{p}[theo]{Proposition}
\newtheorem{lem}[theo]{Lemma}
\newtheorem{de}[theo]{Definition}
\newtheorem{co}[theo]{Corollary}

\newcommand{\n}[1]{\ensuremath{\left\| #1 \right\|}}
\newcommand{\ck}{\ensuremath{\mathcal{K}_E(K)}}
\newcommand{\ch}{\ensuremath{\mathcal{K}_E(H)}}
\newcommand{\ci}{\ensuremath{\mathcal{K}\left(l^2(I)\right)}}
\newcommand{\cx}[1]{\ensuremath{\mathcal#1}}
\newcommand{\cc}[1]{\ensuremath{\mathcal{C}(#1)}}
\newcommand{\cd}[2]{\ensuremath{\mathcal{C}\left(#1,#2\right)}}
\newcommand{\h}[1]{Theorem \ref{#1}}
\newcommand{\pr}[1]{Proposition \ref{#1}}
\newcommand{\prr}{[C] Proposition 5.6.}
\newcommand{\cor}[1]{Corollary {\ref{#1}}}
\newcommand{\lm}[1]{Lemma{ \ref{#1}}}
\newcommand{\dd}[1]{Definition {\ref{#1}}}
\newcommand{\mae}[5]{\ensuremath{#1:#2\longrightarrow #3\,, \quad #4 \longmapsto  #5}}
\newcommand{\mad}[4]{\ensuremath{#1\longrightarrow #2, \quad #3 \longmapsto   #4}}
\newcommand{\mac}[3]{\ensuremath{#1:#2\longrightarrow #3}}
\newcommand{\qedd}{\eqno{\rule{3mm}{3mm}}}
\newcommand{\qed}{\hfill {\rule{3mm}{3mm}}}
\newcommand{\z}[1]{\ensuremath{\{#1\}}}
\newcommand{\ri}[1]{Hilbert right $#1$-module}
\newcommand{\me}[2]{\ensuremath{\left\{\, #1\;|\; #2 \, \right\}}}
\newcommand{\mea}[2]{\ensuremath{\left\{\, #1\;\rule[-5mm]{0.1mm}{1.2cm}\; #2 \, \right\}}}
\newcommand{\s}[2]{\ensuremath{\left \langle  \, #1\left|\, #2 \right. \, \right \rangle}}
\newcommand{\sd}[2]{\ensuremath{#1\left \langle  \, \cdot \left|\, #2 \right. \, \right \rangle}}
\newcommand{\sa}[2]{\ensuremath{\left \langle \, #1 \,, \, #2 \, \right \rangle}}
\newcommand{\si}[1]{\ensuremath{\sum\limits_{#1}}}
\newcommand{\sii}[2]{\ensuremath{\sum\limits_{#1}^{#2}}}
\newcommand{\pro}[1]{\ensuremath{\prod\limits_{#1}}}
\newcommand{\cb}[1]{\ensuremath{\mathop{\bigcirc \hspace{-2.7mm} |}\limits_{ \hspace{2mm} #1}}}
\newcommand{\cw}[1]{\ensuremath{\mathop{\bigcirc \hspace{-2.7mm} |}\limits_{ \hspace{2mm}#1}^{\;\;W}}}
\newcommand{\la}[1]{\ensuremath{\mathcal L\left(#1\right)}}
\newcommand{\lb}[2]{\ensuremath{\mathcal L_{#1}\left(#2\right)}}
\newcommand{\lbb}[2]{\ensuremath{\mathcal L^{#1}\left(#2\right)}}

\newcommand{\lc}[3]{\ensuremath{\mathcal L_{#1}^{#2}(#3)}}
\newcommand{\lca}[1]{\ensuremath{\mathcal L_E^{#1}(H)}}
\newcommand{\lf}[2]{\ensuremath{\mathcal L({#1},#2)}}
\newcommand{\bk}{\ensuremath{\mathrm{I\! K}}}
\newcommand{\bn}{\ensuremath{\mathrm{I\! N}}}
\newcommand{\bz}{\ensuremath{\mathrm{Z\!\!\!Z}}}
\newcommand{\brs}{\ensuremath{\mathrm{I\! R_+}}}
\newcommand{\br}{\ensuremath{\mathrm{I\! R}}}
\newcommand{\bc}{\ensuremath{\mathrm{I\!\!\! C}}}
\newcommand{\aaaa}[4]{\ensuremath{ \left\{ \begin{array}{c@{\quad \mbox{if} \quad}c} #1&#2 \\
#3&#4 \end{array}} \right.}
\renewcommand{\labelenumi}{\alph{enumi})}
\renewcommand{\labelenumii}{$\alph{enumi}_{\arabic{enumii}}$)}

\begin{document}
\begin{center}
\Large {\bfseries
Compact operators\\ on Hilbert right modules}
\end{center}
\begin{flushright}
\large \bfseries
Corneliu Constantinescu
\end{flushright}
\vspace{1cm}

\begin{abstract}
We generalize some results on compact operators on Hilbert spaces to "compact" operators on some Hilbert right W*-modules. We present in this frame the Schatten decomposition of the compact operators, the trace, the Banach $\cx{L}^p$-spaces and their duality, the Hilbert-Schmitt operators, and the integral operators as an example of Hilbert-Schmitt operators. 
\end{abstract}

AMS classification code: 46L08
 
Keywords: compact operators, \ri{W^*}s
\begin{center}
\addtocounter{section}{-1}
\section{Notation and terminology}
\end{center}

In general we use the notation and terminology of [C]. 
In the sequel we give a list of such notation and terminology from [C] used in this paper.
\renewcommand{\labelenumi}{\arabic{enumi}.}
\begin{enumerate}
\item $\bk$ denotes the field of real numbers $\br$ or the field of complex numbers $\bc$. The whole theory is developed in parallel for the real and complex case, but the proofs coincide. $\bz$ denotes the set of integers, $\bn$ denotes the set of natural numbers ($0\not\in \bn$) and we put for every $n\in \bn$,
$$\bn_n:=\me{k\in \bn}{k\leq n}.$$
An initial segment of $\bn$ is a subset $N$ of $\bn$ such that given $m\in \bn$ and $n\in N$, with $m<n$, then $m\in N$. $\brs$ denotes the set of positive real numbers ($0\in \brs$).
\item If $A$ is a set then $id_A$ denotes the identity map of $A$.
\item If $E$ is a Banach space then $E^{\#}$ denotes the unit ball of $E$:
$$ E^{\#}:= \me{x\in E}{\|  x\|  \leq 1}.$$
If $T$ is a compact space then $\cd{T}{E}$ denotes the Banach space of continuous maps $T\rightarrow E$ (endowed with the supremum norm). We put $\cc{T}:=\cd{T}{\bk}$. 
\item Let $E$ be a C*-algebra. We denote by $E_+$ the set of positive elements of $E$ and put $E_+^{\#}:=E_+\cap E^{\#}$. If $E$ is unital then $1_E$ denotes its unit. For $x\in E$, $\sigma (x)$ denotes the spectrum of $x$.
\item If $I$ is a set, then $l^2(I)$ denotes the Hilbert space of square summable families in $\bk$ indexed by $I$, $\la{l^2(I)}$ the W*-algebra of operators 
$$l^2(I)\rightarrow l^2(I),$$ 
and $\ci$ the C*-subalgebra of $\la{l^2(I)}$ of compact operators. 
\item $\delta _{ij}$ denotes Kronecker's symbol:
$$\delta _{ij}:= \aaaa{1}{i=j}{0}{i\not=j}.$$
\item Let $E$ be a C*-algebra and $H$ a \ri{E}. We denote by $\la{H}$ the Banach space of operators $H\rightarrow H$, by $\lb{E}{H}$ its Banach subspace of adjointable operators, which is a C*-algebra, and by $\ch$ the C*-subalgebra of $\lb{E}{H}$ of "compact" operators. For all $\xi ,\eta \in H$ we denote by $\s{\xi }{\eta }$ their scalar product and put
$$\mae{\xi \s{\cdot}{\eta }}{H}{H}{\zeta }{\xi \s{\zeta }{\eta }}.$$
\end{enumerate}

\begin{center}
\fbox{\parbox{12cm}{Throughout this paper we denote by $T$ a compact hyperstonian space ([C] Definition 1.7.2.12), by $E:=\cx{C}(T)$ the C*-algebra of continuous scalar valued functions on $T$ (by [C] Theorem 4.4.4.22 c$\Rightarrow$a, $E$ is a W*-algebra), by $K$ a selfdual Hilbert right $E$-module, by $(p_\iota )_{\iota \in I}$ a family of orthogonal projection of $E$ such that $K$ is isomorphic to $\cw{\iota \in I}p_\iota E$ ([C] Proposition 5.6.4.10 a)), and put $H:=\cb{\iota \in I}p_\iota E$ (by [C] Proposition 5.6.4.1 c), $H$ is a Hilbert right $E$-module)}} 
\end{center}

\begin{center}
\section{The C*-algebra \ch}
\end{center}

\renewcommand{\labelenumi}{\alph{enumi})}

\begin{de}\label{12}
We define \boldmath $\psi$ \unboldmath and for every $t \in T $, \boldmath $\psi_t$ \unboldmath and \boldmath $\varphi_t$ \unboldmath by
$$\mae{\psi }{l^2(I)}{H}{\zeta }{(\zeta _\iota p_\iota )_{\iota\in I}}$$
$$\mae{\psi_t }{H}{l^2(I)}{\xi }{(\xi _\iota(t ))_{\iota\in I}},$$
$$\mae{\varphi _t }{\lb{E}{H}}{\la{l^2(I)}}{u}{\psi _t \circ u\circ \psi }.$$
\end{de}

\begin{p}\label{48}
For every $\xi \in H$ the map
$$\mad{T}{l^2(I)}{t }{\psi _t \xi }$$
is continuous.
\end{p} 

Let $\varepsilon >0$. There is a finite subset $J$ of $I$ such that
$$\si{\iota\in I\setminus J}|\xi _\iota(t )|^2<\varepsilon $$
for all $t \in T $. For $t ,t '\in T $,
$$\n{\psi _t \xi -\psi _{t '}\xi }^2=\si{\iota\in I}|\xi _\iota(t )-\xi_\iota(t ') |^2\leq $$
$$\leq \si{\iota\in J}| \xi_\iota(t ) -\xi_\iota(t ') |^2+2\si{\iota\in I\setminus J}|\xi_\iota(t ) |^2+2\si{\iota\in I\setminus J}|\xi_\iota(t ') |^2\leq $$
$$\leq \si{\iota\in J}|\xi_\iota(t ) -\xi_\iota(t ') |^2+4\varepsilon ,$$
and this implies the assertion.\qed

\begin{p}\label{45} 
Let $t \in T $. 
\begin{enumerate}
\item $\psi _t\circ \psi \circ \psi _t=\psi _t$ .
\item For $\xi ,\eta \in H $ and $\zeta \in l^2(I)$,
$$\s{\psi _t \xi }{\psi _t \eta }=(\s{\xi }{\eta })(t ),$$
$$\s{\psi _t\xi }{\zeta }=\s{\psi _t\xi }{\psi _t\psi \zeta }=(\s{\xi }{\psi \zeta })(t).$$
\item For every $u\in \lb{E}{H}$,
$$\psi _t \circ u\circ \psi \circ \psi _t =\psi _t \circ u.$$
\item For $u,v\in \lb{E}{H}$,
$$\varphi _t (uv)=(\varphi _t u)(\varphi _t v).$$
\item For every $u\in \lb{E}{H}$,
$$\varphi _t u^*=(\varphi _t u)^*.$$
\item For $\xi ,\eta \in H$,
$$\varphi _t (\xi \s{\cdot }{\eta })=(\psi _t \xi) \s{\cdot }{\psi _t \eta }.$$
\end{enumerate}
\end{p}

a) and b) are easy to see.

c) For $\xi \in H$, by a), $\psi _t (\xi -\psi \psi _t \xi )=0$.
Let $\varepsilon >0$. By \pr{48}$\!$, there is a neighborhood $U$ of $t $ such that $\n{\psi _{t '}(\xi -\psi \psi _{t} \xi )}<\varepsilon $ for every $t '\in U$. Let $x\in E_+^{\#}$ with $x(t )=1$ and $x=0$ on $T \setminus U$. Then $\n{(\xi -\psi \psi _t \xi )x}<\varepsilon $ and
$$\n{(u(\xi -\psi \psi _t \xi ))x}=\n{u((\xi -\psi \psi _t \xi )x)}\leq \varepsilon \n{u},$$
$$\n{\psi _t (u(\xi -\psi \psi _t \xi ))}=\n{\psi _t ((u(\xi -\psi \psi _t \xi ))x)}\leq \varepsilon \n{u}.$$
Since $\varepsilon $ is arbitrary,
$$\psi _t u\xi =\psi _t u\psi \psi _t \xi \,,\qquad \psi _t \circ u=\psi _t \circ u\circ \psi \circ \psi _t .$$

d) For $\zeta \in l^2(I)$, by c),
$$(\varphi _t u)(\varphi _t v)\zeta =\psi _t u\psi \psi _t v\psi \zeta =\psi _t uv\psi \xi =(\varphi _t (uv))\zeta ,$$
$$(\varphi _t u)(\varphi _t v)=\varphi _t (uv).$$

e) For $\xi ,\eta \in l^2(I)$, by b),
$$\s{\xi }{(\varphi _t u)^*\eta }=\s{(\varphi _t u)\xi }{\eta }=\s{\psi _t u\psi \xi }{\psi _t \psi \eta }=(\s{u\psi \xi }{\psi \eta })(t )=$$
$$=(\s{\psi \xi }{u^*\psi \eta })(t )=\s{\psi _t \psi \xi }{\psi _t u^*\psi \eta }=\s{\xi }{(\varphi _t u^*)\eta },$$
$$(\varphi _t u)^*=\varphi _t u^*.$$

f) For $\zeta \in l^2(I)$, by b),
$$\varphi _t (\xi \s{\cdot }{\eta })\zeta =\psi _t ((\xi \s{\cdot }{\eta })\psi \zeta )=\psi _t (\xi \s{\psi \zeta }{\eta })=(\psi _t \xi )(\s{\psi \zeta }{\eta })(t )=$$
$$=(\psi _t \xi )\s{\psi _t \psi \zeta }{\psi _t \eta }=(\psi _t \xi )\s{\zeta }{\psi _t \eta }=((\psi _t \xi )\s{\cdot }{\psi _t \eta })\zeta ,$$
$$\varphi _t (\xi \s{\cdot }{\eta })=(\psi _t \xi )\s{\cdot }{\psi _t \eta }.\qedd$$

\begin{co}\label{19.1}
$\rule{0mm}{0mm}$
\begin{enumerate}
\item The map
$$\mad{\lb{E}{H}}{\pro{t \in T }\la{l^2(I)}}{u}{(\varphi _t u)_{t \in T }}$$
is an injective C*-homomorphism.
\item $u\in \lb{E}{H}$ is positive iff $\varphi _tu$ is positive for all $t\in T$.
\end{enumerate}
\end{co}

a) By \pr{45} d),e), the map
$$\mad{\lb{E}{H}}{\pro{t \in T }\la{l^2(I)}}{u}{(\varphi _t u)_{t \in T }}$$
is a C*-homomorphism. Let $u\in \lb{E}{H}$ such that $\varphi _t u=0$ for all $t \in T $. For $\xi \in H$ and $t \in T $, by \pr{45} c),
$$\psi _t u\xi =\psi _tu\psi \psi _t\xi =(\varphi _t u)\psi _t \xi =0\,,\qquad u\xi =0\,,\qquad u=0,$$
so the above map is injective.

b) follows from a).\qed

\begin{p}\label{26}
$\rule{0mm}{0mm}$\begin{enumerate}
\item For every $u\in \ch$ the map
$$\mae{\bar{u} }{T }{\ci}{t }{\varphi _t u}$$
is continuous.
\item The map
$$\mad{\ch}{\cd{T}{\ci}}{u}{\bar{u} }$$
is an injective C*-homomorphism.
\end{enumerate}
\end{p}

a) Let $\xi ,\eta \in H$ and $t ,t '\in T $. By \pr{45} f),
$$\varphi _t (\xi \s{\cdot }{\eta })-\varphi _{t'} (\xi \s{\cdot }{\eta })=(\psi _t \xi )\s{\cdot }{\psi _t \eta }-(\psi _{t'} \xi )\s{\cdot }{\psi _{t'} \eta }=$$
$$=(\psi _t \xi )\s{\cdot }{\psi _t \eta }-(\psi _{t} \xi )\s{\cdot }{\psi _{t'} \eta }+(\psi _t \xi )\s{\cdot }{\psi _{t'} \eta }-(\psi _{t'} \xi )\s{\cdot }{\psi _{t'} \eta }=$$
$$=(\psi _t \xi )\s{\cdot }{\psi _t \eta -\psi _{t '}\eta }+(\psi _t \xi -\psi _{t '}\xi )\s{\cdot }{\psi _{t '}\eta },$$
so by \prr 5.2 a),
$$\n{\varphi _t (\xi \s{\cdot }{\eta })-\varphi _{t'} (\xi \s{\cdot }{\eta })}\leq$$
$$\leq  \n{(\psi _t \xi )\s{\cdot }{\psi _t \eta -\psi _{t '}\eta }}+\n{(\psi _t \xi -\psi _{t '}\xi )\s{\cdot }{\psi _{t '}\eta }}\leq $$
$$\leq \n{\psi _t \xi }\n{\psi _t \eta -\psi _{t '}\eta }+\n{\psi _t \xi -\psi _{t '}\xi }\n{\psi _{t '}\eta }\leq $$
$$\leq \n{\xi }\n{\psi _t \eta -\psi _{t '}\eta }+\n{\psi _t \xi -\psi _{t '}\xi }\n{\eta }.$$
Thus by \pr{48}, the map
$$\mad{T }{\ci}{t }{\varphi _{t }(\xi \s{\cdot }{\eta })}$$
is continuous.

The assertion follows now from the definition of $\ch$ ([C] Definition 5.6.5.3).   

b) follows from a) and \cor{19.1} a).\qed 

\begin{center}
\section{The C*-algebra \cd{T}{\ci}}
\end{center}
\renewcommand{\labelenumi}{\alph{enumi})}

\begin{p}\label{23.1}
Let $u\in \cd{T}{\ci}$ and $n\in \bn$.
\begin{enumerate}
\item The map \boldmath $\theta _n(u)$ \unboldmath defined by 
$$\mae{\theta _n(u)}{T}{\brs}{t}{\theta _n(u(t))}$$
\emph{(with the notation of [C] Definition 6.1.2.1)} is continuous.
\item $\theta _n(u)=\theta _n(u^*)=\theta _n(|u|)$.
\item If $u$ is positive and $f$ is a continuous increasing function on $\brs$ with $f(0)=0$ then $\theta _n(f(u))=f(\theta _n(u))$.
\end{enumerate}
\end{p}

a) follows from [C] Corollary 6.1.2.8.

b) follows from [C] Theorem 6.1.3.1 b).

c) follows from [C] Corollary 6.1.2.16.\qed

\begin{p}\label{10.1}
If $\xi ,\eta \in H$ then
$$\mae{\theta _1(\sd{\xi }{\eta })}{T}{\brs}{t}{\n{\psi _t\xi }\n{\psi _t\eta }}$$
and $\theta _n(\sd{\xi }{\eta })=0$ for all $n\in \bn\setminus \z{1}$.
\end{p}

For $n\in \bn$ and $t\in T$, by \pr{45} f), \pr{26} a), and [C] Proposition 6.1.2.3,
$$(\theta _n(\sd{\xi }{\eta }))(t)=\theta _n(\varphi _t(\sd{\xi }{\eta }))=$$
$$=\theta _n(\sd{(\psi _t\xi )}{\psi _t\eta })=\aaaa{\n{\psi _t\xi }\n{\psi _t\eta }}{n=1}{0}{n\not=1}.\qedd$$

\begin{de}\label{23.1a}
We put for every $\xi \in K$ and $t\in T$,
\begin{center}
\boldmath ${{\xi (t)}}$\unboldmath $:=(\xi _\iota (t))_{\iota \in I}\in l^2(I).$
\end{center}
We put for every $u\in \cd{T}{\ci}$ and $n\in \bn$
$${\bf{U_n(u)}}:=\me{t\in T}{\theta _n(u(t))\not=0},$$
$$\mae{{\bf{e_n(u)}}}{T}{\bk}{t}{\aaaa{1}{t\in \overline{U_n(u)}}{0}{t\in T\setminus \overline{U_n(u)}}}.$$
A sequence $(\xi _n)_{n\in \bn}$ in $K$ is called \boldmath $u\!\!$ \unboldmath {\bf{-orthonormal}} if for all $m,n\in \bn$, $m\leq n$,
$$\s{\xi _m}{\xi _n}=\delta _{m,n}e_n(u)$$
and the map
$$\mad{U_n(u)}{l^2(I)}{t}{\xi _n(t)}$$
is continuous.
We extend the above notation and terminology to $u\in \ch$ by using \emph{\pr{26} a)}.
\end{de}

If $\xi \in H$ then $\xi (t)=\psi _t\xi $ for all $t\in T$.

\begin{p}\label{28.1}
Let $u\in \cd{T}{\ci}$ and let $(\xi _n)_{n\in \bn}$ be a $u$-orthonormal sequence in $K$.
\begin{enumerate}
\item $U_n(u)$ is the union of a sequence of pairwise disjoint clopen sets of $T$ for every $n\in \bn$.
\item $\sd{\xi _n}{\xi _n}$ is an orthogonal projection of $\ck$ for every $n\in \bn$ and
$$(\sd{\xi _m}{\xi _m})(\sd{\xi _n}{\xi _n})=0$$
for all distinct $m,n\in \bn$.
\end{enumerate}
\end{p}

a) If we denote for every $k\in \bz$ by $U_k$ the closure of the interior of the set
$$\me{t\in T}{2^k\leq \theta _n(u(t))<2^{k+1}}$$
then $(U_k)_{k\in \bz}$ is a countable set of pairwise disjoint clopen sets of $T$ the union of which is $T$.

b) For all $m,n\in \bn$, $m\leq n$,
$$(\sd{\xi _m}{\xi _m})(\sd{\xi _n}{\xi _n})=\sd{(\xi _m\s{\xi _n}{\xi _m})}{\xi _n}=\delta _{m,n}\;\sd{\xi _m}{\xi _n}.\qedd$$

\begin{p}\label{14}
Let $u$ be a selfadjoint element of $\cd{T }{\ci}$.
\begin{enumerate}
\item For every $t \in T $ there is a representation
$$u(t )=\si{\alpha \in \sigma (u(t ))}\alpha \pi _{t ,\alpha },$$
where for every $\alpha \in \sigma (u(t ))$, $\pi _{t ,\alpha }$ is the orthogonal projection of $l^2(I)$ onto $Ker\,(\alpha 1-u(t ))$ \emph{(here $1=id_{l^2(I)}$)}  and $\pi _{t ,\alpha }\pi _{t ,\beta }=0$ for all distinct $\alpha ,\beta \in \sigma (u(t ))$. 
\item Let $t \in T $, $\alpha \in \sigma (u(t ))$, $\alpha \not=0$, $\varepsilon >0$, and $U$ a neighborhood of $\alpha $ such that $\sigma (u(t ))\cap \bar{U}=\z{\alpha } $ and $|\alpha -\beta| \leq \frac{|\alpha |\varepsilon }{2}$ for all $\beta \in U$. Then there is a neighborhood $V$ of $t $ such that for every $t '\in V$,
$$\n{\si{\beta \in \sigma (u(t '))\cap U}\beta \pi _{t ',\beta }-\alpha \pi _{t ,\alpha }}<\varepsilon\,, \qquad 
\n{\si{\beta \in \sigma (u(t '))\cap U}\pi _{t ',\beta }-\pi _{t ,\alpha }}<\varepsilon .$$
\end{enumerate}
\end{p}

a) follows from [C] Theorem 5.5.6.1 a$\Rightarrow$c\&e.

b) Let $U'$ be a neighborhood of $\sigma (u(t ))\setminus \z{\alpha }$ such that $\bar{U}\cap \bar{U'}=\emptyset $. By [C] Corollary 2.2.5.2, there is a neighborhood $W$ of $t $ such that $\sigma (u(t '))\subset U\cup U'$ for all $t '\in W$. Let $f\in \cx{C}(\bk)_+$, $0\leq f\leq 1$, $f=1$ on $\bar{U} $, and $f=0$ on $\bar{U'} $. By [C] Proposition 4.1.3.20, the map
$$\mad{T }{\ci}{t' }{f(u(t '))}$$
is continuous. Thus there is a neighborhood $V$ of $t $, $V\subset W$, such that for every $t '\in V$,
$$\n{f(u(t '))-f(u(t ))}<\inf\left\{\varepsilon ,\frac{|\alpha |\varepsilon }{2}\right\}.$$
By [C] Theorem 5.5.6.1 a$\Rightarrow $f,
$$f(u(t ))=\alpha \pi _{t ,\alpha }\,,\qquad f(u(t '))=\si{\beta \in \sigma (u(t '))\cap U}\beta \pi _{t ',\beta }.$$
It follows
$$\n{\si{\beta \in \sigma (u(t '))\cap U}\beta \pi _{t ',\beta }-\alpha \pi _{t ,\alpha }}=\n{f(u(t '))-f( u(t ))}<\inf\left\{\varepsilon ,\frac{|\alpha |\varepsilon }{2}\right\} ,$$
$$\n{\si{\beta \in \sigma (u(t '))\cap U}\pi _{t ',\beta }-\pi _{t ,\alpha }}=\frac{1}{|\alpha |}\n{\si{\beta \in \sigma (u(t '))\cap U}\alpha \pi _{t ',\beta }-\alpha \pi _{t ,\alpha }}\leq $$
$$\leq \frac{1}{|\alpha |}\n{\si{\beta \in \sigma (u(t '))\cap U}(\alpha -\beta )\pi _{t ',\beta }}+\frac{1}{|\alpha |}\n{\si{\beta \in \sigma (u(t '))\cap V}\beta \pi _{t ',\beta }-\alpha \pi _{t ,\alpha }}\leq $$
$$\leq \frac{|\alpha -\beta |}{|\alpha |}+\frac{1}{|\alpha |}\frac{|\alpha |\varepsilon }{2}\leq \varepsilon \,.\qedd$$

\begin{lem}\label{14b}
Let $\mac{\eta }{T }{l^2(I)}$ be a map such that the map
$$\mad{T }{\ci}{t }{\eta (t )\s{\cdot }{\eta (t )}}$$
is continuous. Let $t _0\in T $ with $\eta (t _0)\not=0$ and put
$$U:=\me{t \in T }{\s{\eta (t _0)}{\eta (t )}\not=0},$$
$$\mae{\xi }{U}{l^2(I)}{t }{\frac{\s{\eta (t _0)}{\eta (t )}}{|\s{\eta (t _0)}{\eta (t )}|}}\eta (t ).$$
Then $U$ is an open neighborhood of $t_0$, $\xi $ is continuous, $\xi (t _0)=\eta (t _0)$, and
$$\xi (t )\s{\cdot }{\xi (t )}=\eta (t )\s{\cdot }{\eta (t )}$$
for all $t \in U$. 
\end{lem}

The map
$$\mad{T }{\brs}{t }{\s{\eta (t )\s{\eta (t _0)}{\eta (t )}}{\eta (t _0)}}=|\s{\eta (t )}{\eta (t _0)}|^2$$
is continuous so
$$\lim_{t \rightarrow t _0}|\s{\eta (t _0)}{\eta (t )}|=|\s{\eta (t _0)}{\eta (t _0)}|\not=0.$$
Thus $U$ is an open neighborhood of $t_0$, $\xi $ is continuous, $\xi (t _0)=\eta (t _0)$, and
$$\xi (t )\s{\cdot }{\xi (t )}=\eta (t )\s{\cdot }{\eta (t )}$$
for all $t \in U$.\qed

\begin{co}\label{14a}
Let $u$ be a positive element of $\cd{T }{\ci}$.
\begin{enumerate}
\item For every $t \in T $ there are an initial segment $N_t$ of $\bn$ and an orthonormal family $(\eta _{t,n})_{n\in N_t}$ in $l^2(I)$ such that $\eta _{t,n}=0$ for all $t\in T\setminus U_n(u)$ and 
$$u(t )=\si{n\in N_t}\theta _n(u(t ))\;\sd{\eta _{t,n}}{\eta _{t,n}}.$$
\item Let $t_0\in T$ such that $N_{t_0}$ is finite and let $U$ be a neighborhood of $t_0$ such that $N_t=N_{t_0}$ for all $t\in U$. Then there is a neighborhood $V$ of $t _0$ and for every $n\in N_{t_0}$ a continuous map
$$\mac{\xi_n }{V}{l^2(I)}$$
such that for every $t\in V$, $(\xi _n(t))_{n\in N_{t_0}}$ is an orthonormal family in $l^2(I)$ and
$$\sd{\xi _n(t)}{\xi _n(t)}=\sd{\eta _{t,n}}{\eta _{t,n}}.$$ 
\end{enumerate}
\end{co}

a) follows from [C] Corollary 6.1.2.13 a$\Rightarrow$ b\&c.

b) follows\pr{14} b) and \lm{14b}.\qed

\begin{p}\label{25a}
If $u$ is a positive element of $\cd{T}{\ci}$ then there is a $u$-orthonormal sequence $(\xi _n)_{n\in \bn}$ in $K$ such that for every 
$$t\in T\setminus \bigcup_{n\in \bn}\left(\overline{U_n(u)}\setminus U_n(u)\right), $$
$$u(t)=\si{n\in \bn}\theta _n(u(t))\,\sd{(\xi _n(t))}{\xi _n(t)}\qquad \emph{(in\; \ci)}.$$
\end{p}
 
By \cor{14a} a), for every $t\in T$ there is an initial segment $N_t$ of $\bn$ and an orthonormal family $(\xi _{t,n})_{n\in N_t}$ in $l^2(I)$ such that $\xi _{t,n}=0$ for all $t\in T\setminus \overline{U_n(u)}$ and $n\in N_t$ and
$$u(t)=\si{n\in N_t}\theta _n(u(t))\;\sd{\xi _{t,n}}{\xi _{t,n}} \qquad \mbox{(in\; \ci)}.$$
For every $k\in \bn$, let $f_k\in \cc{\brs}$ with $0\leq f_k\leq 1$, $f_k=0$ on $[0,\frac{1}{2k}]$, $f_k=1$ on $[\frac{1}{k},\infty ]$. By [C] Proposition 4.1.3.20, for every $k\in \bn$ the map
$$\mad{T}{\ci}{t}{f_k(u(t))}$$
is continuous. By \pr{23.1} c), for $t\in T$,
$$f_k(u(t))=\si{n\in N_t}f_k(\theta _n(u(t)))\;\sd{\xi _{t,n}}{\xi _{t,n}}.$$
By \pr{23.1} a), $(\theta _n(u))_{n\in \bn}$ is a decreasing sequence of continuous real functions on $T$ with infimum $0$, so by Dini's theorem it converges uniformly to $0$ on $T$. Thus by \pr{23.1} c), for every $k\in \bn$ there is an $m\in \bn$ such that
$$\theta _m(f_k(u))=0.$$
Since $T$ is hyperstonian and since $U_n(u)$ is the union of a sequence of clopen sets of $T$ (\pr{28.1} a)), we may assume (by \cor{14a} b)) that for every $n\in \bn$ there is a $\xi _n\in K$ such that the map
$$\mad{U_n(u)}{l^2(I)}{t}{\xi _n(t)}$$
is continuous, with $\s{\xi _n}{\xi _n}=e_n(u)$ and $\sd{\xi _n(t)}{\xi _n(t)}=\sd{\xi _{t,n}}{\xi _{t,n}}$ for all $t\in T$. Moreover for $m,n\in \bn$, $m<n$, and $t\in U_n(u)$,
$$\sd{\xi _m(t)}{\xi _n(t)}\s{\xi _n(t)}{\xi _m(t)}=(\sd{\xi _m(t)}{\xi _m(t)})\circ (\sd{\xi _n(t)}{\xi _n(t)})=$$
$$=(\sd{\xi _{t,m}}{\xi _{t,m}})\circ (\sd{\xi _{t,n}}{\xi _{t,n}})=\sd{\xi _{t,m}}{\xi _{t,n}}\s{\xi _{t,n}}{\xi _{t,m}}=0.$$
By \pr{10.1}, $\s{\xi _n(t)}{\xi _m(t)}=0$ so $\s{\xi_ n}{\xi _m}=0$. Thus $(\xi _n)_{n\in \bn}$ is $u$-orthonormal.\qed

\begin{theo}\label{28.1a}
Let $u\in \ch\;(\subset \ck)$.
\begin{enumerate}
\item If $u$ is positive then there is a $u$-orthonormal sequence $(\xi _n)_{n\in \bn}$ in $K$ such that
$$u=\si{n\in \bn}\theta _n(u)\;\sd{\xi _n}{\xi _n}\qquad \emph{(in\;\ck)}.$$
In this case $u\xi _n=\theta _n(u)\xi _n\in H$ for all $n\in \bn$.
\item There are $u$-orthonormal sequences $(\xi _n)_{n\in \bn}$ and $(\eta _n)_{n\in \bn}$ in $K$ such that
$$u=\si{n\in \bn}\theta _n(u)\;\sd{\xi _n}{\eta _n}\qquad \emph{(in\;\ck)}.$$
The above identities are called {\bf{Schatten decomposition of}} \boldmath{$u$}\unboldmath.
\end{enumerate}
\end{theo}

By [C] Theorem 5.6.3.5 b), $\lb{E}{K}$ is a W*-algebra with $\stackrel{...}{K}$ as predual.

a) Let $(\xi _n)_{n\in \bn}$ be the $u$-orthonormal sequence in $K$ defined in \pr{25a}. By \pr{28.1} b), for $k,m\in \bn$, $k\leq m$,
$$\sii{n=k}{m}\theta _n(u)\;\sd{\xi _n}{\xi _n}\leq \theta _k(u)\sii{n=k}{m}\sd{\xi _n}{\xi _n}\leq \theta _k(u),$$
so the sequence $(\theta _n(u)\;\sd{\xi _n}{\xi _n})_{n\in \bn}$ is summable in $\ck$. By \pr{25a} (and [C] Definition 5.6.3.2),
$$u=\si{n\in \bn}\theta _n(u)\;\sd{\xi _n}{\xi _n}$$
in $\lb{E}{K}$ with respect to its weak topology associated to the duality 
$$\sa{\lb{E}{K}}{\stackrel{...}{K}},$$
so
$$u=\si{n\in \bn}\theta _n(u)\;\sd{\xi _n}{\xi _n}\qquad \mbox{(in \ck)}.$$

From $u\xi _n=\theta _n(u)\xi _n$ it follows 
$$(u\xi _n)(t)=\theta _n(u(t))\xi _n(t)$$
for all $t\in T$. Thus the map
$$\mad{T}{\bk}{t}{\s{(u\xi _n)(t)}{(u\xi _n)(t)}=\theta _n(u(t))^2\s{\xi _n(t)}{\xi _n(t)}}$$
is continuous and $u\xi _n\in H$.

b) By a) (and \pr{23.1} b)), there is a $u$-orthonormal sequence $(\eta _n)_{n\in \bn}$ in $K$ such that
$$|u|=\si{n\in \bn}\theta _n(u)\;\sd{\eta _n}{\eta _n}\qquad \mbox{(in\;\ck)}.$$
Let $u=w|u|$ be the polar representation of $u$ ([C] Theorem 4.4.3.1). Then
$$u=\si{n\in \bn}\theta _n(u)\;\sd{(w\eta _n)}{\eta _n}\qquad \mbox{(in\;\ck)}.$$
For $m,n\in \bn$, $m\leq n$, since $w^*w$ is the carrier of $|u|$ and 
$$|u|\eta _n=\theta _n(u)\eta _n,$$
$$\theta _n(u)\s{w\eta _n}{w\eta _n}=\s{\eta _n}{w^*w\theta _n(u)\eta _n}=\s{\eta _n}{w^*w|u|\eta _n}=$$
$$=\s{\eta _n}{|u|\eta _n}=\theta _n(u)\s{\eta _n}{\eta _n},$$
so by \pr{28.1} b),
$$\s{w\eta _m}{w\eta _n}=\delta _{m,n}e_n(u).$$
Thus if we put $\xi _n:=w\eta _n$ for every $n\in \bn$ then 
$$u=\si{n\in \bn}\theta _n(u)\;\sd{\xi _n}{\eta _n}\qquad \mbox{(in\;\ck)}.$$

Let $n\in \bn$. Since the map
$$\mad{U_n(u)}{l^2(I)}{t}{\eta _n(t)}$$
is continuous, the map
$$\mad{U_n(u)}{l^2(I)}{t}{u\eta _n(t)}$$
is also continuous. From 
$$u\eta _n=\theta _n(u)\xi _n,$$
it follows that the map
$$\mad{U_n(u)}{l^2(I)}{t}{\xi _n(t)}$$
os continuous. Thus $(\xi _n)_{n\in \bn}$ is a $u$-orthonormal sequence in $K$.\qed

\begin{p}\label{6.2}
Let $A$ be a dense set of $T$ and $(\theta _n)_{n\in \bn}$ be a decreasing sequence in $E_+$ such that
$$\lim_{n\rightarrow \infty }\theta _n(t)=0$$
for every $t\in A$. Let further $(\xi _{n,t})_{(n,t)\in \bn\times A}$ and $(\eta _{n,t})_{(n,t)\in \bn\times A}$ be families in $l^2(I)$ such that
 $(\xi _{n,t})_{n\in N_t}$ and $(\eta _{n,t})_{n\in N_t}$ are orthonormal families in $l^2(I)$ for all $t\in A$, where
$$N_t:=\me{n\in \bn}{\xi _{n,t}\not=0}=\me{n\in \bn}{\eta _{n,t}\not=0}.$$
If for an $u\in \ch$,
$$\varphi _tu=\si{n\in \bn}\theta _n(t)\,\sd{\xi _{n,t}}{\eta _{n,t}}\qquad \left(\emph{in\; \ci}\right)$$
for all $t\in A$ then $\theta _n(u)=\theta _n$ for all $n\in \bn$.
\end{p}

By [C] Proposition 6.1.2.11, for $t\in A$,
$$(\theta _n(u))(t)=\theta _n(\varphi _tu)=\theta _n(t),$$
so $\theta _n(u)=\theta _n$, since $\theta _n(u)$ is continuous (\pr{23.1} a)).\qed

\begin{co}\label{2.2}
Let $u\in \ch$ and let
$$u=\si{n\in \bn}\theta _n(u)\,\sd{\xi _n}{\eta _n}$$
be a Schatten decomposition of $u$.
\begin{enumerate}
\item $$u^*=\si{n\in \bn}\theta _n(u)\,\sd{\eta_n}{\xi_n}$$
is a Schatten decomposition of $u^*$.
\item $\theta _n(u^*u)=\theta _n(u)^2$ for every $n\in \bn$ and
$$u^*u =\si{n\in \bn}\theta _n(u)^2\,\sd{\eta_n}{\eta_n}$$
is a Schatten decomposition of $u^*u$.
\item Let $N$ be a subset of $\bn$ and
$$v:=\si{n\in N}\theta _n(u)\,\sd{\xi _n}{\eta _n}.$$
If $M$ is an initial segment of $\bn$ and $\mac{f}{M}{N}$ is an increasing bijective map then
$$\theta _n(v)=\aaaa{\theta _{f(n)}(u)}{n\in M}{0}{n\in \bn\setminus M}.$$ 
\end{enumerate}
\end{co}

a) By [C] Proposition 5.6.5.2 a),
$$u^*=\si{n\in \bn}\theta _n(u)\,\sd{\eta_n}{\xi_n}\qquad \mbox{(in\;\ck )}$$
and the assertion follows from \pr{23.1} b).

b) By a), for $n\in \bn$,
$$u^*\xi _n=\si{m\in \bn}\theta _m(u)\eta _m\s{\xi _n}{\xi _m}=\theta _n(u)\eta _n,$$
so
$$u^*u=\si{n\in \bn}\theta _n(u)\,\sd{(u^*\xi _n)}{\eta _n}=\si{n\in \bn}\theta _n(u)^2\,\sd{\eta _n}{\eta _n}.$$
If we put
$$\mae{\eta '_n}{T}{l^2(I)}{t}{\aaaa{\eta _n(t)}{t\in U_n(u)}{0}{T\setminus U_n(u)}}$$
for every $n\in \bn$ then
$$\varphi _t(u^*u)=\si{n\in \bn}(\theta _n(u)^2)(t)\,\sd{\eta '_n(t)}{\eta' _n(t)}$$
for all $t\in T$ and the assertion follows from \pr{6.2}.

c) The above defined sequence $(\theta _n(v))_{n\in\bn}$ is decreasing and converges to $0$. Put
$$A:=T\setminus \bigcup_{n\in\bn}\left(\overline{U_n(u)}\setminus U_n(u)\right) $$
and for every $n\in \bn$ and $t\in A$,
$$\xi _{n,t}:=\aaaa{\xi _{f(n)}(t)}{n\in M}{0}{n\in \bn\setminus M},\quad \eta _{n,t}:=\aaaa{\eta _{f(n)}(t)}{n\in M}{0}{n\in \bn\setminus M}.$$
Then for $t\in A$,
$$\varphi _t(v)=\si{n\in N}(\theta _n(u))(t)\;\sd{\xi _n(t)}{\eta _n(t)})=$$
$$=\si{n\in M}(\theta _{f(n)}(u))(t)\;\sd{\xi _{f(n)}(t)}{\eta _{f(n)}(t)}=
\si{n\in \bn}(\theta _n(v))(t)\;\sd{\xi _{n,t}}{\eta _{n,t}}$$
and the assertion follows from \pr{6.2}.\qed

\begin{center}
\section{The Banach spaces $\lc{E}{p}{H}$}
\end{center}

\begin{de}\label{25}
We denote for every $p\in [1,\infty [$ by \boldmath{$\lc{E}{p}{H}$}\unboldmath\, the set of $u\in \ch$ for which the sequence $(\theta _n^p)_{n\in \bn}$ is summable in $E$ and define \boldmath $\n{\cdot }_p$\unboldmath\, by
$$\mae{\n{\cdot }_p}{\lc{E}{p}{H}}{\brs}{u}{\n{\si{n\in \bn}\theta _n(u)^p}}^{\frac{1}{p}}.$$
 Moreover we put \boldmath{$\lc{E}{\infty }{H}$}\unboldmath \,$:=\lb{E}{H}$, \boldmath{$\lc{E}{0}{H}$}\unboldmath \,$:=\ch$, and define \boldmath $\n{\cdot }_0$\unboldmath\, by
$$\mae{\n{\cdot }_0}{\lc{E}{0}{H}}{\brs}{u}{\n{u}=\n{\theta _1(u)}}.$$
\end{de}

\begin{p}\label{16.1}
Let $u,v\in \ch$, $0\leq u\leq v$.
\begin{enumerate}
\item $\theta _n(u)\leq \theta _n(v)$ for all $n\in \bn$.
\item If $p,q\in [1,\infty [$, $p\leq q$, and $v\in \lc{E}{p}{H}$ then $u\in \lc{E}{q}{H}$.
\end{enumerate} 
\end{p}

a) By \cor{19.1} b), for $t\in T$, $0\leq \varphi _tu\leq \varphi _tv$ and this implies $\theta _n(\varphi _tu)\leq \theta _n(\varphi _tv)$ ([C] Definition 6.1.2.1).

b) Let $\zeta \in H$. By [C] Theorem 5.6.1.11 c),
$$\s{v\zeta }{\zeta }^q=\s{v\zeta }{\zeta }^{q-p}\s{v\zeta }{\zeta }^p\leq \n{v}^{q-p}\n{\zeta }^{2(q-p)}\s{v\zeta }{\zeta }^p,$$
so $\theta _n(v)^q\leq \n{v}^{q-p}\theta _n(v)^p$ for all $n\in \bn$ ([C] Definition 6.1.2.1) and therefore $v\in \lc{E}{q}{H}$. By a), $u\in \lc{E}{q}{H}$.\qed

\begin{p}\label{4.1}
Let $p\in [1,\infty [$.
\begin{enumerate}
\item If $u\in \ch_+$ then
$$u\in \lc{E}{p}{H}\Longleftrightarrow u^p\in \lc{E}{1}{H}\Longrightarrow \n{u}_p^p=\n{u^p}_1.$$
\item If $u\in \ch$ then
 $$u\in \lc{E}{p}{H}\Longleftrightarrow u^*\in \lc{E}{p}{H}\Longleftrightarrow |u|\in \lc{E}{p}{H}\Longrightarrow $$
$$  \Longrightarrow \n{u}_p=\n{u^*}_p=\n{\,|u|\,}_p.$$
\end{enumerate}
\end{p}

a) By \pr{23.1} c), $\theta _n(u^p)=\theta _n(u)^p$ for all $n\in \bn$.

b) follows from \pr{23.1} b).\qed

\begin{de}\label{31.1}
We denote by \boldmath{$\Omega$}\unboldmath\, the set of sequences $(\zeta _n)_{n\in \bn}$ in $K$ such that:
\renewcommand{\labelenumi}{\arabic{enumi}.}
\begin{enumerate}
\item For every $n\in \bn$ there is a closed nowhere dense set $F_n$ of $T$ such that the map
$$\mad{T\setminus F_n}{l^2(I)}{t}{\zeta _n(t)}$$
is continuous.
\item $(\zeta _n(t))_{n\in N_t}$ is an orthonormal family in $l^2(I)$ for all $t\in T$, where 
$$N_t:=\me{n\in \bn}{\zeta _n(t)\not=0}.$$
\end{enumerate}
\end{de}
\renewcommand{\labelenumii}{$\alph{enumi}_{\arabic{enumii}}$)}

\begin{p}\label{31}
Let  $p\in [1,\infty [$.
\begin{enumerate}
\item If $u\in \lc{E}{p}{H}$ then
$$\si{n\in \bn}\theta _n(u)^p=\sup\mea{\si{n\in \bn}|\s{u\zeta _n}{\zeta '_n}|^p}{(\zeta _n)_{n\in \bn} ,(\zeta' _n)_{n\in \bn}\in \Omega }.$$
\item If $u$ is a positive element of $\lc{E}{p}{H}$ then
$$\si{n\in \bn}\theta _n(u)^p=\sup\mea{\si{n\in \bn}\s{u\zeta _n}{\zeta _n}^p}{(\zeta _n)_{n\in \bn} \in \Omega }.$$
\end{enumerate}
 \end{p}

a) Let 
$$u=\si{n\in N}\theta _n(u)\,\sd{\xi _n}{\eta _n}$$
be a Schatten decomposition of $u$ and put for every $n\in \bn$
$$\mae{\xi '_n}{T}{l^2(I)}{t}{\aaaa{\xi _n(t)}{t\in U_n(u)}{0}{t\in T\setminus U_n(u)}},$$
$$\mae{\eta '_n}{T}{l^2(I)}{t}{\aaaa{\eta _n(t)}{t\in U_n(u)}{0}{t\in T\setminus U_n(u)}}.$$
Then $(\xi' _n)_{n\in \bn},(\eta' _n)_{n\in \bn}\in \Omega $, so
$$\si{n\in \bn}\theta _n(u)^p=\si{n\in \bn}|\s{u\eta _n}{\xi _n}|^p=\si{n\in \bn}|\s{u\eta' _n}{\xi' _n}|^p\leq$$
$$\leq \sup\mea{\si{\lambda \in L}|\s{u\zeta _\lambda }{\zeta '_\lambda }|^p}{(\zeta _\lambda)_{\lambda \in L} ,(\zeta' _\lambda )_{\lambda \in L}\in \Omega }\,. $$

Let $(\zeta _n)_{n\in \bn} ,(\zeta '_n)_{n\in \bn}\in \Omega $ and $t\in T$. We put for all $m,n\in \bn$,
$$\alpha _{m,n}:=\s{\xi _n(t)}{\zeta '_m(t)}\s{\zeta _m(t)}{\eta _n(t)}.$$
If $m\in \bn$ then
$$\si{n\in \bn}|\alpha _{m,n}|=\si{n\in \bn}|\s{\xi _n(t)}{\zeta '_m(t)}\s{\zeta _m(t)}{\eta _n(t)}|\leq $$
$$\leq \left(\si{n\in \bn}|\s{\xi _n(t)}{\zeta '_m(t)}|^{\,2}\right)^{\frac{1}{2}}\left(\si{n\in \bn}|\s{\zeta _m(t)}{\eta _n(t)}|^{\,2}\right)^{\frac{1}{2}}\leq$$
$$\leq  \n{\zeta '_m(t)}\n{\zeta _m(t)}\leq 1.$$
If $n\in \bn$ then
$$\si{m\in \bn}|\alpha _{m,n}|=\si{m\in \bn}|\s{\xi _n(t)}{\zeta '_m(t)}\s{\zeta _m(t)}{\eta _n(t)}|\leq $$
$$\leq \left(\si{m\in \bn}|\s{\xi _n(t)}{\zeta '_m(t)}|^{\,2}\right)^{\frac{1}{2}}\left(\si{m\in \bn}|\s{\zeta _m(t)}{\eta _n(t)}|^{\,2}\right)^{\frac{1}{2}}\leq$$
$$\leq  \n{\xi _n(t)}\n{\eta _n(t)}\leq 1.$$
For $m\in \bn$,
$$\s{(\varphi _tu)\zeta _m(t)}{\zeta '_m(t)}=\si{n\in \bn}\theta _n(\varphi _t(u))\s{\xi _n(t)}{\zeta '_m(t)}\s{\zeta _m(t)}{\eta _n(t)}.$$
By [C] Lemma 6.1.3.9,
$$\si{n\in \bn}|\s{(\varphi _tu)\zeta _n(t)}{\zeta '_n(t)}|^p\leq \si{n\in \bn}\theta _n(\varphi _tu)^p.$$
Since
$$\s{(\varphi _tu)\zeta _n(t)}{\zeta '_n(t)}=(\s{u\zeta _n}{\zeta '_n})(t)$$
for all $t\in T\setminus \bigcup\limits _{n\in \bn}F_n$, we get
$$\si{n\in \bn}|\s{u\zeta _n}{\zeta '_n}|^p\leq \si{n\in \bn}\theta _n(u)^p,$$
$$\sup\mea{\si{n\in \bn}|\s{u\zeta _n}{\zeta '_n}|^p}{(\zeta _n)_{n\in \bn} ,(\zeta' _n)_{n\in \bn}\in \Omega }\leq \si{n\in \bn}\theta _n(u)^p.$$

b) The proof is similar to the proof of a).\qed

\begin{theo}\label{31a} 
Let $p\in [1,\infty [$.
\begin{enumerate}
\item $\lc{E}{p}{H}$ is a vector subspace of $\ch$ and the map
$$\mad{\lc{E}{p}{H}}{\brs}{u}{\n{u}_p}$$
is a norm. We always consider $\lc{E}{p}{H}$ endowed with this norm.
\item $\lc{E}{p}{H}$ is complete.
\item If $u\in \lca{p}$ and
$$u=\si{n\in \bn}\theta _n(u)\,\sd{\xi _n}{\eta _n}$$
is a Schatten decomposition of $u$ with $\xi _n,\eta _n\in H$ for all $n\in \bn$ then the above sum converges in $\lca{p}$.
\end{enumerate}
\end{theo}

a) Let $u,v\in \lc{E}{p}{H}$. By [C] Proposition 6.1.2.5, for $n\in \bn$,
$$\theta _{2n-1}(u+v)\leq \theta _n(u)+\theta _n(v),$$
$$\theta _{2n}(u+v)\leq \theta _n(u)+\theta _{n+1}(v),$$
so
$$\theta _{2n-1}(u+v)^p\leq (\theta _n(u)+\theta _n(v))^p\leq 2^{p-1}(\theta _n(u)^p+\theta _n(v)^p),$$
$$\theta _{2n}(u+v)^p\leq (\theta _n(u)+\theta _{n+1}(v))^p\leq 2^{p-1}(\theta _n(u)^p+\theta _{n+1}(v)^p).$$
Thus $u+v\in \lc{E}{p}{H}$. Let $(\xi _n )_{n \in \bn}\,,(\eta _n )_{n \in \bn}\in \Omega $. By \pr{31} a),
$$\left(\si{n \in \bn}|\s{(u+v)\xi _n }{\eta _n }|^p\right)^{\frac{1}{p}}=\left(\si{n \in \bn}|\s{u\xi _n }{\eta _n }+\s{v\xi _n }{\eta _n }|^p\right)^{\frac{1}{p}}\leq $$
$$\leq \left(\si{n \in \bn}|\s{u\xi _n }{\eta _n }|^p\right)^{\frac{1}{p}}+\left(\si{n \in \bn}|\s{v\xi _n }{\eta _n }|^p\right)^{\frac{1}{p}}\leq \n{u}_p+\n{v}_p,$$
$$\n{u+v}_p\leq \n{u}_p+\n{v}_p.$$
 
b) Let $(u_n)_{n\in \bn}$ be a Cauchy sequence in $\lc{E}{p}{H}$. Then $(u_n)_{n\in \bn}$ converges in $\ch$ to a $u$. Let $\varepsilon >0$. There is an $n_0\in \bn$ such that 
$$\n{u_m-u_n}_p<\varepsilon $$
for all $m,n\in \bn\setminus \bn_{n_0}$. Let $(\xi _k )_{k \in \bn}\,,(\eta _k )_{k \in \bn}\in\Omega $.  By a) and \pr{31} a),
$$\n{\si{k \in \bn}|\s{(u_m-u_n)\xi _k }{\eta _k }|^p}\leq \n{u_m-u_n}_p^p<\varepsilon^p $$
for all $m,n\in \bn\setminus \bn_{n_0}$. Hence
$$\n{\si{k \in \bn}|\s{(u_n-u)\xi _k }{\eta _k }|^p}<\varepsilon ^p$$
for all $n\in \bn\setminus \bn_{n_0}$. By a) and \pr{31} a), again,
$$u_n-u\in \lc{E}{p}{H}\,,\qquad u\in \lc{E}{p}{H}\,,\qquad\n{u_n-u}_p<\varepsilon $$
for all $n\in \bn\setminus \bn_{n_0}$. Thus $(u_n)_{n\in \bn}$ converges to $u\in \lc{E}{p}{H}$ and $\lc{E}{p}{H}$ is complete.

c) By \cor{2.2} c), for $n_0\in \bn$,
$$\n{\sii{n=n_0}{\infty }\theta _n(u)\,\sd{\xi _n}{\eta _n}}_p=\left(\sii{n=n_0}{\infty }\theta _n(u)^p\right)^{\frac{1}{p}}.\qedd$$

\begin{co}\label{1.1}
If $p\in [1,\infty [$, $u\in \lc{E}{p}{H}$, and $v,w\in \lb{E}{H}$ then
$$vuw\in \lc{E}{p}{H}\,,\qquad\qquad \n{vuw}_p\leq \n{v}\n{u}_p\n{w}.$$
\end{co}

By \pr{45} d) and [C] Corollary 6.1.3.13 a), for $t\in T$ and $n\in \bn$,
$$\theta _n(\varphi _t(vuw))=\theta _n((\varphi _tv)(\varphi _tu)(\varphi _tw))\leq$$
$$\leq  \n{\varphi _tv}\theta _n(\varphi _tu)\n{\varphi _tw}\leq \n{v}\theta _n(\varphi _tu)\n{w}$$
and the assertion follows.\qed 

\begin{co}\label{1.1a}
Let $p\in \z{0}\cup [1,\infty [$ and let $q\in [1,\infty ]$ be the conjugate exponent of $p$.
\begin{enumerate}
\item If $u\in \lc{E}{p}{H}$ and $v\in \lc{E}{q}{H}$ then
$$uv,vu\in \lc{E}{1}{H},$$
$$\n{uv}_1\leq \n{u}_p\n{v}_q,\qquad \n{vu}_1\leq \n{u}_p\n{v}_q\quad \mbox{\bf{(\emph{H\"older\; inequality})}}.$$ 
\item For every $u\in \lca{p}$ there is a $v\in \lca{q}$ such that
$$\n{uv}_1=\n{vu}_1=\n{u}_p\n{v}_q.$$
\end{enumerate}
\end{co}

a) By \cor{1.1} we may assume $p\in ]1,\infty [$. By [C] Corollary 6.1.2.7, for $n\in \bn$,
$$\theta _{2n-1}(uv)\leq \theta _n(u)\theta _n(v),\qquad\qquad \theta _{2n}(uv)\leq \theta _n(u)\theta _{n+1}(v),$$
so for $N\subset \bn$,
$$\si{n\in N}\theta _{2n-1}(uv)\leq \si{n\in N}\theta _n(u)\theta _n(v)\leq \left(\si{n\in N}\theta _n(u)^p\right)^{\frac{1}{p}}\left(\si{n\in N}\theta _n(v)^q\right)^{\frac{1}{q}},$$
$$\si{n\in N}\theta _{2n}(uv)\leq \si{n\in N}\theta _n(u)\theta _{n+1}(v)\leq \left(\si{n\in N}\theta _n(u)^p\right)^{\frac{1}{p}}\left(\si{n\in N}\theta _{n+1}(v)^q\right)^{\frac{1}{q}}.$$
Thus $(\theta _n(uv))_{n\in \bn}$ is summable in $E$ and $uv\in \lc{E}{1}{H}$. By [C] Theorem 6.1.3.21, for $t\in T$,
$$\si{n\in \bn}\theta _n(\varphi _t(uv))\leq \left(\si{n\in \bn}\theta _n(\varphi _t(u))^p\right)^{\frac{1}{p}}\left(\si{n\in \bn}\theta _n(\varphi _t(v))^q\right)^{\frac{1}{q}},$$
$$\n{uv}_1\leq \n{u}_p\n{v}_q.$$
The assertion for $vu$ follows.

b) Let 
$$u=\si{n\in\bn}\theta _n(u)\,\sd{\xi _n}{\eta _n}$$
be a Schatten decomposition of $u$. If $p=1$ then we may take $v=id_H$. Assume $p=0$. Put
$$v:=\sd{\eta _1}{\xi _1}.$$
By \pr{10.1}, $v\in \lca{1}$, $\n{v}_1=1$,
$$uv=\si{n\in\bn}\theta _n(u)\,(\sd{\xi _n}{\eta _n})(\sd{\eta _1}{\xi _1})=$$
$$=\si{n\in\bn}\theta _n(u)\,\sd{\xi _n\s{\eta _1}{\eta _n}}{\xi _1}=\theta _1(u)\,\sd{\xi _1}{\xi _1},$$ 
$$vu=\si{n\in\bn}\theta _n(u)\,(\sd{\eta _1}{\xi _1})(\sd{\xi _n}{\eta _n})=$$
$$=\si{n\in\bn}\theta _n(u)\,\sd{\eta _1\s{\xi _n}{\xi _1}}{\eta _n}=\theta _1(u)\,\sd{\eta _1}{\eta _1}.$$
Thus (by \pr{10.1})
$$\n{uv}_1=\n{vu}_1=\n{\theta _1(u)}=\n{u}_p\n{v}_q.$$

Assume now $p\in ]1,\infty [$. Put
$$v:=\si{n\in\bn}\theta _n(u)^{\frac{p}{q}}\,\sd{\eta _n}{\xi _n}\qquad \mbox{(in\;\ck)}.$$
By \cor{2.2} c), $\theta _n(v)=\theta _n(u)^{\frac{p}{q}}$ for every $n\in\bn$ so
$$v\in \lca{q},\qquad\qquad \n{v}_q^q=\n{u}_p^p.$$
For $n\in\bn$,
$$u\eta _n=\theta _n(u)\xi _n,\qquad v\xi _n=\theta _n(u)^{\frac{p}{q}}\eta _n,$$
so
$$uv=\si{n\in\bn}\theta _n(u)^{\frac{p}{q}+1}\,\sd{\xi _n}{\xi _n},\qquad vu=\si{n\in\bn}\theta _n(u)^{1+\frac{p}{q}}\,\sd{\eta _n}{\eta _n}.$$
By \cor{2.2} c),
$$\theta _n(uv)=\theta _n(vu)=\theta _n(u)^{\frac{p}{q}+1}=\theta _n(u)^p,$$
$$\n{uv}_1=\n{vu}_1=\si{n\in\bn}\theta _n(u)^p=\n{u}_p^p=$$
$$=\n{u}_p\n{u}_p^{p-1}=\n{u}_p\n{v}_q^{\frac{q}{p}(p-1)}=\n{u}_p\n{v}_q.\qedd$$

\begin{center}
\section{The trace}
\end{center}

\begin{p}\label{1.2a}
Let $(\theta _n)_{n\in \bn}$ be a summable sequence in $E_+$ and let $(\xi _n)_{n\in \bn}$ and $(\eta _n)_{n\in \bn}$ be sequences in $K^{\#}$. 
\begin{enumerate}
\item $(\theta _n\,\sd{\xi _n}{\eta _n})_{n\in \bn}$ is summable in $\ck$; we put
$$u:=\si{n\in \bn}\theta _n\,\sd{\xi _n}{\eta _n}.$$
\item For every Fourier basis $A$ of $K$ \emph{([C] Definition 5.6.3.11)}
$$\si{n\in \bn}\theta _n\s{\xi _n}{\eta _n}=\si{\zeta \in A}\s{u\zeta }{\zeta }.$$
\end{enumerate}
\end{p}

a) By [C] Proposition 5.6.5.2 a),
$$\n{\sd{\xi _n}{\eta _n}}\leq \n{\xi _n}\n{\eta _n}\leq 1$$
for every $n\in \bn$.

b) For $\zeta \in A$,
$$\s{u\zeta }{\zeta }=\si{n\in \bn}\theta _n\s{\xi _n}{\zeta }\s{\zeta }{\eta _n}.$$
By [C] Theorem 5.6.3.13 f), since the above sum converges uniformly,
$$\si{\zeta \in A}\s{u\zeta }{\zeta }=\si{\zeta \in A}\,\si{n\in \bn}\theta _n\s{\xi _n}{\zeta }\s{\zeta }{\eta _n}=$$
$$=\si{n\in \bn}\theta _n\si{\zeta \in A}\s{\xi _n}{\zeta }\s{\zeta }{\eta _n}=\si{n\in \bn}\theta _n\s{\xi _n}{\eta _n}.\qedd$$

\begin{de}\label{1.2}
Let $u\in \lc{E}{1}{H}$ and let 
$$u:=\si{n\in \bn}\theta _n(u)\;\sd{\xi _n}{\eta _n}$$
be a Schatten decomposition of $u$. We put
$${\bf{tr\,u}}:=\si{n\in \bn}\theta _n(u)\s{\xi _n}{\eta _n}\in E$$
and call it {\bf{the trace of}} \boldmath{$u$}\unboldmath \; (by \emph{\pr{1.2a} b)} the trace of $u$ does not depend on the chosen Schatten decomposition of $u$).
\end{de}

\begin{co}\label{3.1a}
Given $u\in \lb{E}{K}$ and $\xi ,\xi ',\eta ,\eta '\in K$,
$$tr\,(\sd{\xi }{\eta })=\s{\xi }{\eta },$$
$$tr\,(u\circ (\sd{\xi }{\eta }))=\s{u\xi }{\eta }=tr\,((\sd{\xi }{\eta })\circ u),$$
$$tr\,((\sd{\xi }{\eta })\circ (\sd{\xi '}{\eta '}))=\s{\xi }{\eta '}\s{\xi '}{\eta }.$$
\end{co}

[C] Proposition 5.6.5.2 d), e).\qed

\begin{p}\label{4.2}
We put for all $u\in \lb{E}{H}$ and $x\in E$,
$$\mae{{\bf{ux}}}{H}{H}{\xi }{(u\xi )x=u(\xi x)}.$$
Then $ux\in \lb{E}{H}$, $(ux)^*=u^*x^*$, and $\n{ux}\leq \n{u}\n{x}$ for all $u\in \lb{E}{H}$ and $x\in E$,
\end{p}

For $\xi ,\eta \in H$,
$$\s{(ux)\xi }{\eta }=\s{(u\xi )x}{\eta }=\s{u\xi }{\eta }x=$$
$$=\s{\xi }{u^*\eta }x=
\s{\xi }{(u^*\eta )x^*}=\s{\xi }{(u^*x^*)\eta },$$
so $ux\in \lb{E}{H}$ and $(ux)^*=u^*x^*$. For $\xi \in H$,
$$\n{(ux)\xi }=\n{(u\xi )x}\leq \n{u\xi }\n{x}\leq \n{u}\n{\xi }\n{x},$$
so $\n{ux}\leq \n{u}\n{x}$.\qed

\begin{co}\label{3.1c}  
The map
$$\mad{\lc{E}{1}{H}}{E}{u}{tr\,u}$$
is linear, involutive, positive, and continuous with norm 1 \emph{(\h{31a} a))} and
$$\n{tr\,u}=\n{u}_1$$
for every positive element of $\lc{E}{1}{H}$. Moreover for all $u\in \lc{E}{1}{H}$ and $x\in E$ \emph{(\pr{4.2})},
$$tr\,(ux)=(tr\,u)x.$$
\end{co}

tr is linear (\pr{1.2a} b)), involutive (\cor{2.2} a)), and continuous with norm at most $1$ ([C] proposition 5.6.5.2 a)). By \dd{1.2}, tr is positive and 
$$\n{tr\,u}=\n{u}_1$$
If $A$ is a Fourier basis of $K$ then by \pr{1.2a} b) ,
$$tr\,(ux)=\si{\zeta \in A}\s{(ux)\zeta }{\zeta }=\left(\si{\zeta \in A}\s{u\zeta }{\zeta }\right)x=(tr\,u)x.\qedd$$

\begin{co}\label{4.1a}
If $u\in \ch_+$ and $p\in [1,\infty [$ then
$$u\in \lc{E}{p}{H}\Longleftrightarrow u^p\in \lc{E}{1}{H}\Longrightarrow \n{u}_p=(tr\,u^p)^{\frac{1}{p}}.$$
\end{co}

By \pr{4.1} a), $u\in \lc{E}{p}{H}$ iff $u^p\in \lc{E}{1}{H}$ and 
$$\n{u}_p^p=\n{u^p}_1.$$
By \cor{3.1c},
$$\n{u}_p=(tr\,u^p)^{\frac{1}{p}}\,.\qedd$$

\begin{p}\label{4.1b}
If $u\in \lc{E}{1}{H}$ and $v\in \lb{E}{H}$ then \emph{(\cor{1.1})}
$$tr\,(uv)=tr\,(vu).$$
\end{p}

Let
$$u=\si{n\in \bn}\theta _n(u)\,\sd{\xi _n}{\eta _n}$$
be a Schatten decomposition of $u$. By [C] Proposition 5.6.5.2 d),e) (and [C] Theorem 5.6.4.7 d)),
$$tr\,(vu)=tr\,\si{n\in \bn}\theta _n(u)\,\sd{(v\xi _n)}{\eta _n}=\si{n\in \bn}\theta _n(u)\,\s{v\xi _n}{\eta _n}=$$
$$=\si{n\in \bn}\theta _n(u)\,\s{\xi _n}{v^*\eta _n}
=tr\,\si{n\in \bn}\theta _n(u)\,\sd{\xi _n}{v^*\eta _n}=tr\,(uv).\qedd$$

\begin{center}
\section{Hilbert-Schmidt operators}
\end{center}

\begin{de}\label{1.1b}
The elements of $\lc{E}{2}{H}$ are called {\bf{Hilbert-Schmidt operators on}} \boldmath $H$. \unboldmath
\end{de}

\begin{p}\label{2.1}
$\lc{E}{2}{H}$ endowed with the exterior multiplications \emph{(\cor{1.1})}
$$\mad{\lb{E}{H}\times \lc{E}{2}{H}}{\lc{E}{2}{H}}{(w,u)}{wu},$$
$$\mad{\lc{E}{2}{H}\times \lb{E}{H}}{\lc{E}{2}{H}}{(u,w)}{uw}$$
and with the inner-product \emph{(\cor{1.1a} a))}
$$\mae{\s{\cdot }{\cdot }}{\lc{E}{2}{H}\times \lc{E}{2}{H}}{\lb{E}{H}}{(u,v)}{v^*u}$$
is a unital Hilbert $\lb{E}{H}$-module \emph{([C] Definition 5.6.1.4)}.
\end{p}

For $u,v\in \lc{E}{2}{H}$ and $w\in \lb{E}{H}$,
$$\s{u}{v}^*=(v^*u)^*=u^*v=\s{v}{u},$$
$$\s{uw}{v}=v^*(uw)=(v^*u)w=\s{u}{v}w,$$
$$\s{wu}{v}=v^*(wu)=(w^*v)^*u=\s{u}{w^*v},$$
$$\s{wu}{wu}=u^*w^*wu\leq \n{w}^2u^*u=\n{w}^2\s{u}{u},$$
$$1_{\lb{E}{H}}u=u.$$
Moreover if $\bk=\br$,
$$(\s{u}{u}+\s{v}{v},\s{v}{u}-\s{u}{v})=(u^*u+v^*v,u^*v-v^*u)=(u,v)^*(u,v)$$
is a positive element of the complexification of $\lb{E}{H}$.\qed

\begin{p}\label{1.1c}
For every $u\in \ch$, 
$$u\in \lc{E}{2}{H}\Longleftrightarrow u^*u\in \lc{E}{1}{H}\Longrightarrow \n{u^*u}_1=\n{u}_2^2.$$ 
\end{p}

If $u\in \lc{E}{2}{H}$ then by \cor{2.2} b), $u^*u\in \lc{E}{1}{H}$ and 
$$\n{u^*u}_1=\si{n\in \bn}\theta _n(u^*u)=\si{n\in \bn}\theta _n(u)^2=\n{u}_2^2.$$
If $u^*u\in \lc{E}{1}{H}$ then by \cor{2.2} b), $(\theta _n(u)^2)_{n\in \bn}$ is summable in $E$ so $u\in \lc{E}{2}{H}$.\qed

\begin{theo}\label{4.1d}
$\rule{0mm}{0mm}$
\begin{enumerate}
\item $u,v\in \lc{E}{2}{H}\Longrightarrow v^*u\in \lc{E}{1}{H}$.
\item $\lc{E}{2}{H}$ endowed with the exterior multiplication \emph{(\pr{4.2})}
$$\mad{\lc{E}{2}{H}\times E}{\lc{E}{2}{H}}{(u,x)}{ux}$$
 and with the inner-product \emph{(a))}
$$\mae{\s{\cdot }{\cdot }}{\lc{E}{2}{H}\times \lc{E}{2}{H}}{E}{(u,v)}{tr\,(v^*u)}$$
is a Hilbert right $E$-module with norm $\n{\cdot }_2$.
\item $u,v\in \lc{E}{2}{H}\Longrightarrow \s{u}{v}=\s{v^*}{u^*}$.
\end{enumerate}
\end{theo}

a) follows from the H\"older inequality.

b) For $u,v\in \lc{E}{2}{H}$ and $x\in E$, by \cor{3.1c} and \pr{1.1c},
$$\s{ux}{v}=tr\,(v^*ux)=tr\,(v^*u)x=\s{u}{v}x,$$
$$\s{u}{v}=tr\,(v^*u)=(tr\,(u^*v))^*=\s{v}{u}^*,$$
$$\s{u}{u}=tr\,(u^*u)\in E_+,\qquad \n{\s{u}{u}}=\n{u}_2^2.$$

c) By \pr{4.1b},
$$\s{u}{v}=tr\,(v^*u)=tr\,(uv^*)=\s{v^*}{u^*}.\qedd$$

\begin{center}
\section{Duals of $\lc{E}{p}{H}$-spaces}
\end{center}

\begin{p}\label{8.2}
Let $p\in [1,\infty [$ and let $\cx{F}$ be the set of $u\in \lc{E}{p}{H}$ for which there is a Schatten decomposition
$$u=\si{n\in \bn}\theta _n(u)\,\sd{\xi _n}{\eta _n}$$
such that $(\xi _n)_{n\in \bn}$ and $(\eta _n)_{n\in \bn}$ are sequences in $H$. Then $\cx{F}$ is dense in $\lc{E}{p}{H}$.
\end{p}

Let $u\in \lc{E}{p}{H}$ and let
$$u=\si{n\in \bn}\theta _n(u)\,\sd{\xi _n}{\eta _n}$$
be a Schatten decomposition of $u$. We put for all $n,k\in \bn$,
$$U_{n,k}:=\mea{t\in T}{\theta _n(t)>\frac{1}{kn^2}},$$
$$\mae{e_{n,k}}{T}{\bk}{t}{\aaaa{1}{t\in \overline{U_{n,k}}}{0}{t\in T\setminus \overline{U_{n,k}}}},$$
$$u_k:=\si{n\in \bn}\theta _n(u)\,\sd{\xi _n}{\eta _ne_{n,k}}=\si{n\in \bn}(\theta _n(u)e_{n,k})\,\sd{(\xi _ne_{n,k})}{\eta _ne_{n,k}}.$$
For $k\in \bn$,
$$u-u_k=\si{n\in \bn}\theta _n(u)\,\sd{\xi _n}{\eta _n(1_E-e_{n,k})}=$$
$$=\si{n\in \bn}(\theta _n(u)(1_E-e_{n,k}))\,\sd{(\xi _n(1_E-e_{n,k}))}{\eta _n(1_E-e_{n,k})}.$$
By \pr{6.2}, for $n,k\in \bn$,
$$\theta _n(u-u_k)=\theta _n(u)(1_E-e_{n,k})\leq \frac{1}{kn^2},$$
so $(\theta _n(u-u_k)^p)_{n\in \bn}$ is summable in $E$ and
$$\si{n\in \bn}\theta _n(u-u_k)^p\leq \frac{1}{k^p}\si{n\in \bn}\frac{1}{n^{2p}}.$$
Thus $(u_k)_{k\in \bn}$ converges to $u$ 
in $\lc{E}{p}{H}$ and this proves the assertion since $u_k\in \cx{F}$ for every $k\in \bn$.\qed

\begin{theo}\label{14.1b}
Let $p\in \z{0}\cup [1,\infty [$, $q\in [1,\infty ]$ the conjugate exponent of $p$, and $\lf{\lc{E}{p}{H}}{E}$ the involutive Banach space of operators from $\lc{E}{p}{H}$ to $E$ \emph{([C] Proposition 2.3.2.22 a))}, the involution being defined for every $\phi \in \lf{\lc{E}{p}{H}}{E}$ by
$$\mae{\phi ^*}{\lca{p}}{E}{u}{\left(\phi (u^*)\right)}^*.$$
Further let $\cx{G}$ be the set of $\phi \in \lf{\lc{E}{p}{H}}{E}$ such that
\renewcommand{\labelenumi}{\arabic{enumi}.}
\begin{enumerate}
\item $u\in \lca{p},\; x\in E\Longrightarrow \phi (ux)=\phi (u)x$
\item For $\xi \in H$,
$$(\phi (\sd{\xi }{e_\iota }))_{\iota \in I},(\phi^* (\sd{\xi }{e_\iota }))_{\iota \in I}\in H,$$
where for every $\iota \in I$,
$$e_\iota :=(\delta _{\iota ,\lambda }1_E)_{\lambda \in I}\;(\in H).$$
\end{enumerate}
\renewcommand{\labelenumi}{\alph{enumi})}
\begin{enumerate}
\item $\cx{G}$ is an involutive vector subspace of $\lf{\lc{E}{p}{H}}{E}$.
\item If we put for every $v\in \lca{q}$ \emph{(by the H\"older inequality and \pr{4.1b})}
$$\mae{\tilde{v} }{\lca{p}}{E}{u}{tr\,(uv)=tr\,(vu)}$$
then $\tilde{v}\in \cx{G} $ and the map
$$\mae{\Psi }{\lca{q}}{\cx{G}}{v}{\tilde{v} }$$
is an isomorphism of involutive Banach spaces.
\end{enumerate}
\end{theo}

a) is easy to see.

b) For $u\in \lca{p}$, by \cor{3.1c} and the H\"older inequality,
$$\n{\tilde{v}(u) }=\n{tr\,(uv)}\leq \n{uv}_1\leq \n{u}_p\n{v}_q,$$
so $\n{\tilde{v} }\leq \n{v}_q$ and $\tilde{v}\in  \lf{\lc{E}{p}{H}}{E}$. By \cor{3.1c}, for $u\in \lca{p}$ and $x\in E$,
$$\tilde{v}(ux)=tr\,(vux)=tr\,(vu)x=\tilde{v}(u)x.  $$
For $\xi \in H$, by \cor{3.1a},
$$(\tilde{v}(\sd{\xi }{e_\iota }) )_{\iota \in I}=tr\,(v(\sd{\xi }{e_\iota }))_{\iota \in I}=(\s{v\xi }{e_\iota })_{\iota \in I}=v\xi \in H,$$
$$(\widetilde{v^*}(\sd{\xi }{e_\iota }) )_{\iota \in I}=v^*\xi \in H,$$
so $\tilde{v}\in \cx{G} $. $\Psi $ is obviously linear. For $u\in \lca{p}$, by \cor{3.1c},
$$\widetilde{v^*}(u)=tr\,(uv^*)=(tr\,(vu^*))^*=(\tilde{v}(u^*) )^*=\tilde{v}^*(u),  $$
so $\widetilde{v^*}=\tilde{v}^*  $ and $\Psi $ is involutive. Moreover by \cor{1.1a}, $\Psi $ is norm preserving. The only thing we have still to prove is the surjectivity of $\Psi $. 

Let $\phi \in \cx{G}$ and put ([C] Proposition 5.6.5.2 a)) 
$$\mae{v}{H}{H}{\xi }{(\phi (\sd{\xi }{e_\iota }))}_{\iota \in I},$$
$$\mae{w}{H}{H}{\xi }{(\phi^* (\sd{\xi }{e_\iota }))}_{\iota \in I}.$$
For $\xi ,\eta \in H$, by $1.$ and [C] Proposition 5.6.5.2 a),c),
$$\s{v\xi }{\eta }=\si{\iota \in I}\s{v\xi }{e_\iota }\eta ^*_\iota =\si{\iota \in I}\phi (\sd{\xi }{e_\iota })\eta ^*_\iota =\phi (\sd{\xi }{\eta }),$$
$$\n{v\xi }^2=\n{\s{v\xi }{v\xi }}=\n{\phi (\sd{\xi }{v\xi })}\leq \n{\phi }\n{\xi }\n{v\xi },$$
$$\n{v\xi }\leq \n{\phi }\n{\xi }, \qquad \n{v}\leq \n{\phi }.$$
For $\iota ,\lambda \in I$, by [C] Proposition 5.6.5.2 a),
$$\s{ve_\lambda }{e_\iota }=\phi (\sd{e_\lambda }{e_\iota })=\phi (\sd{e_\lambda }{e_\iota })^{**}=$$
$$=(\phi ^*(\sd{e_\iota }{e_\lambda }))^*=\s{we_\iota }{e_\lambda }^*=\s{e_\lambda }{we_\iota }.$$
Thus $v\in \lb{E}{H}$ and $v^*=w$. Let $u\in \lca{p}$ and let 
$$u=\si{n\in\bn}\theta  _n(u)\,\sd{\xi _n}{\eta_n}$$
be a Schatten decomposition of $u$ with $(\xi _n)_{n\in\bn}$ and $(\eta _n)_{n\in\bn}$ sequences in $H$. Then by the above and \h{31a} c),
$$\tilde{v}(u)=\si{n\in\bn}\theta _n(u)\,\tilde{v}(\sd{\xi _n}{\eta _n})=\si{n\in\bn}\theta _n(u)\phi (\sd{\xi _n}{\eta _n})=\phi (u).  $$
By \pr{8.2}, $\tilde{v}=\phi  $ and $\Psi $ is surjective.\qed

\begin{center}
\section{Integral operators}
\end{center}

\begin{center}
\fbox{\parbox{12.8cm}{Throughout this section $S$ is a compact space, $\mu $ a positive Radon measure on $S$, $(h_\iota )_{\iota \in I}$ an orthonormal basis of $L^2(\mu )$, $H:=\cb{\iota \in I}E$, and $w\in \cd{S\times S}{E}$. Moreover $\odot $ denotes the algebraic tensor product}} 
\end{center}

\begin{p}\label{38}
The linear map
$$\mad{L^2(\mu )\odot E}{H}{f\otimes x}{(\s{f}{h_\iota }x)_{\iota \in I}}$$
can be extended to an isomorphism $L^2(\mu )\otimes E\longrightarrow H$ \emph{([L] pages 34-35)} of Hilbert right modules.
\end{p}

We denote by $\Phi $ the above map. For $(f,x),(g,y)\in L^2(\mu )\times E$ and $z\in E$,
$$\s{\Phi (f\otimes x)}{\Phi (g\otimes y)}=\s{(\s{f}{h_\iota }x)_{\iota \in I}}{\s{g}{h_\iota }y)_{\iota \in I}}=$$
$$=\si{\iota \in I}y^*\s{h_\iota }{g}\s{f}{h_\iota }x=
y^*\s{f}{g}x=\s{f\otimes x}{g\otimes y},$$
$$\Phi ((f\otimes x)z)=\Phi (f\otimes (xz))=(\s{f}{h_\iota }(xz))_{\iota \in I}=$$
$$=(\s{f}{h_\iota }x)_{\iota \in I}\,z=(\Phi (f\otimes x))z,$$
i.e. $\Phi $ preserves the inner-product and the right multiplication so it can be extended to a linear map 
$$\mac{\Psi }{L^2(\mu )\otimes E}{H}$$
preserving the inner-product and the right multiplication. Moreover
$$\Psi (h_\lambda \otimes z)=(\delta _{\lambda ,\iota }z)_{\iota \in I}$$
for all $\lambda \in I$ and $z\in E$, so $\Psi $ is surjective.\qed

\begin{lem}\label{40}
The vector subspace of $\cd{S\times S}{E}$ generated by maps of the form
$$\mad{S\times S}{E}{(r,s)}{u(r)v(s)},$$
where $u\in \cd{S}{E}$ and $v\in \cd{S}{\bk}$ is dense in $\cd{S\times S}{E}$.
\end{lem}

Let $\varepsilon >0$. There are finite open coverings $(U_j)_{j\in J}, (V_k)_{k\in K}$ of $S$ such that
$$\n{w(r,s)-w(r',s')}<\varepsilon $$
for all $(j,k)\in J\times K$ and $(r,s),(r',s')\in U_j\times V_k$. Take $r_j\in U_j$ and $s_k\in V_k$ for all $j\in J$ and $k\in K$ and let $(f_j)_{j\in J}$ and $(g_k)_{k\in K}$ be partitions of unity subordinate to the coverings $(U_j)_{j\in J}$ and $(V_k)_{k\in K}$ of $S$, respectively. For $r,s\in S$,
$$\n{w(r,s)-\si{(j,k)\in J\times K}f_j(r)g_k(s)w(r_j,s_k)}=$$
$$=\n{\si{(j,k)\in J\times K}f_j(r)g_k(s)(w(r,s)-w(r_j,s_k))}\leq$$
$$\leq  \si{(j,k)\in J\times K}f_j(r)g_k(s)\varepsilon =\varepsilon .$$
If we put 
$$\mae{u_k}{S}{E}{r}{\si{j\in J}f_j(r)w(r_j,s_k)}$$
and $v_k:=g_k$ for all $k\in K$ then for $r,s\in S$,
$$\si{(j,k)\in J\times K}f_j(r)g_k(s)w(r_j,s_k)=\si{k\in K}\left(\si{j\in J}f_j(r)w(r_j,s_k)\right)g_k(s)=$$
$$=\si{k\in K}u_k(r)v_k(s).\qedd$$

\begin{de}\label{39}
A function $\mac{f}{S\times T}{\bk}$ is called {\bf{E-\boldmath $\mu $\unboldmath -integrable }} if $f(s,\cdot )\in E$ and $f(\cdot ,t)\in \cx{L}^1(\mu )$ for all $(s,t)\in S\times T$ and if the map
$$\mad{T}{\bk}{t}{\int f(\cdot ,t)\,\mathrm{d}\mu }$$
is continuous, i.e. it belongs to $E$. We denote this element of $E$ by

\boldmath $$\int g\,\mathrm{d}\mu =\int g(s)\,\mathrm{d}\mu (s),$$\unboldmath
where
$$\mae{g}{S}{E}{s}{f(s,\cdot )}.$$
\end{de}

\begin{lem}\label{41}
For every $f\in L^2(\mu )$ the map
$$\mae{\tilde{f} }{S}{E}{r}{\int w(r,s)f(s)\,\mathrm{d}\mu (s)}$$
is continuous.
\end{lem}

Let $r_0\in S$ and $\varepsilon >0$. There is a neighborhood $U$ of $r_0$ such that
$$\sup_{s\in S}\n{w(r,s)-w(r_0,s)}<\varepsilon $$
for all $r\in U$. Then for $r\in U$,
$$\n{\tilde{f}(r)-\tilde{f}(r_0)  }=\n{\int (w(r,s)-w(r_0,s))f(s)\,\mathrm{d}\mu (s)}\leq \varepsilon \int |f(s)|\,\mathrm{d}\mu (s).\qedd$$

\begin{lem}\label{42}
We use the notation of \emph{\lm{41}}.
\begin{enumerate}
\item The linear map
$$\mad{L^2(\mu )\odot E}{\cd{S}{E}}{f\odot x}{\tilde{f}x }$$ 
is continuous so it can be extended by continuity to an operator 
$$L^2(\mu) \otimes E\longrightarrow \cd{S}{E}.$$
\item The linear map
$$\mad{L^2(\mu )\odot E}{H}{f\odot x}{\tilde{f}x }$$
is continuous so it can be extended by continuity to an operator 
$$\mac{\tilde{w} }{H}{H}.$$
\end{enumerate} 
\end{lem}

a) Let $(f_j)_{j\in J}$ and $(x_j)_{j\in J}$ be finite families in $L^2(\mu )$ and $E$, respectively. For $r\in S$,
$$\rule[-5mm]{0.1mm}{1.2cm}\left(\si{j\in J}\tilde{f}_jx_j \right)(r)\rule[-5mm]{0.1mm}{1.2cm}=
\rule[-5mm]{0.1mm}{1.2cm}\si{j\in J}\int w(r,s)f_j(s)x_j\,\mathrm{d}\mu (s)\rule[-5mm]{0.1mm}{1.2cm}=$$
$$=\rule[-5mm]{0.1mm}{1.2cm}\int w(r,s)\left(\si{j\in J}f_j(s)x_j\,\mathrm{d}\mu (s)\right)\rule[-5mm]{0.1mm}{1.2cm}\leq $$
$$\int |w(r,s)|\,\rule[-5mm]{0.1mm}{1.2cm}\si{j\in J}f_j(s)x_j\rule[-5mm]{0.1mm}{1.2cm}\,\mathrm{d}\mu (s)\leq \n{w}\int \rule[-5mm]{0.1mm}{1.2cm}\si{j\in J}f_j(s)x_j\rule[-5mm]{0.1mm}{1.2cm}\,\mathrm{d}\mu (s),$$
where
$$\n{w}:=\sup_{r,s\in S}\n{w(r,s)}.$$
Thus
$$\rule[-5mm]{0.1mm}{1.2cm}\left(\si{j\in J}\tilde{f}_jx_j \right)(r)\rule[-5mm]{0.1mm}{1.2cm}\leq \n{w}\mu (S)^{\frac{1}{2}}\left(\int \rule[-5mm]{0.1mm}{1.2cm}\si{j\in J}f_j(s)x_j\rule[-5mm]{0.1mm}{1.2cm}^{\;2}\,\mathrm{d}\mu (s)\right)^{\frac{1}{2}}=$$
$$=\n{w}\mu (S)^{\frac{1}{2}}\left(\si{j.k\in J}x_jx_k^*\int f_j(s)\overline{f_k(s)}\,\mathrm{d}\mu (s)\right)^{\frac{1}{2}}=$$
$$=\n{w}\mu (S)^{\frac{1}{2}}\left(\si{j,k\in J}\s{f_j}{f_k}\s{x_j}{x_k}\right)^{\frac{1}{2}}=$$
$$=\n{w}\mu (S)^{\frac{1}{2}}\s{\si{j\in J}(f_j\otimes x_j)}{\si{j\in J}(f_j\otimes x_j)}^{\frac{1}{2}}\leq$$
$$\leq  \n{w}\mu (S)^{\frac{1}{2}}\n{\si{j\in J}(f_j\otimes x_j)}. $$

b) By [W] T3.13,
$$\cd{S}{E}\approx \cd{S}{\bk}\otimes E$$
and by \pr{38}, $L^2(\mu )\otimes E\approx H$. The assertion follows from the continuity of the inclusion $\cd{S}{\bk}\otimes E\subset L^2(\mu )\otimes E$.\qed

\begin{theo}\label{64}
We use the notation of \emph{\lm{42} b)}. $\tilde{w}\in  \lca{2}$ (i.e. $\tilde{w} $ is a Hilbert Schmitt operator on $H$) and $\tilde{w}^*=\widetilde{{w'}}$, where
$$\mae{w'}{S\times S}{E}{(r,s)}{w(s,r)^*}$$
and $\widetilde{w'} $ is defined similarly to $\tilde{w} $.
\end{theo}

\begin{center}
Step 1 $\tilde{w}\in  \lb{E}{H}$ and $\tilde{w}^*=\widetilde{w'} $
\end{center}

For $(f,x),(g,y)\in L^2(\mu )\times E$,
$$\s{\tilde{w}(f\otimes x) }{g\otimes y}=\int y^*g(r)^*\left(\int w(r,s)f(s)x\,\mathrm{d}\mu (s)\right)\,\mathrm{d}\mu (r)=$$
$$=\int f(s)x\left(\int w(r,s)y^*g(r)^*\,\mathrm{d}\mu (r)\right)\,\mathrm{d}\mu (s)=$$
$$=\int f(s)x\left(\int w(r,s)^*g(r)y\,\mathrm{d}\mu (r)\right)^*\,\mathrm{d}\mu (s)=$$
$$=\int f(s)x\left(\widetilde{w'}(g\otimes y) \right)^*(s)\,\mathrm{d}\mu (s)=\s{f\otimes x}{\widetilde{w'}(g\otimes y) }$$
so $\tilde{w}\in  \lb{E}{H}$ and $\tilde{w}^*=\widetilde{w'} $.

\begin{center}
Step 2 $\tilde{w}\in \ch $
\end{center}

By \lm{40}, we may assume that there are $u\in \cd{S}{E}$ and $v\in \cd{S}{\bk}$ with
$$\mae{w}{S\times S}{E}{(r,s)}{u(r)v(s)}.$$
For $(f,x)\in L^2(\mu )\times E$,
$$\tilde{w}(f\otimes x)=\int u\,v(s)f(s)x\,\mathrm{d}\mu (s)=u\s{f}{\bar{v} }\s{x}{1_E}= $$
$$=u\s{f\otimes x}{\bar{v}\otimes 1_E }=(\sd{u}{\bar{v}\otimes 1_E })(f\otimes x),$$
$$\tilde{w}=\sd{u}{\bar{v}\otimes 1_E }\in \ch. $$ 

\begin{center}
Step 3 $\tilde{w}\in \lca{2} $
\end{center}

For $t\in T$,
$$(w(\cdot ,\cdot ))(t)\in \cd{S\times S}{\bk}\subset L^2(\mu \otimes \mu ),$$
so we consider in the sequel $(w(\cdot ,\cdot ))(t)\in L^2(\mu \otimes \mu ).$

Let $t_0\in T$ and $\varepsilon >0$. There is a neighborhood $U$ of $t_0$ such that
$$\sup_{r,s\in S}|(w(r,s))(t)-(w(r,s))(t_0)|<\varepsilon $$
for all $t\in U$. Then
$$\n{(w(\cdot ,\cdot ))(t)-(w(\cdot ,\cdot ))(t_0)}_2^2=$$
$$=\int |(w(r,s))(t)-(w(r,s))(t_0)|^2\,\mathrm{d}(\mu \otimes \mu )(r,s)\leq \varepsilon ^2\mu (S)^2$$
for all $t\in U$. Thus the map
$$\mad{T}{L^2(\mu \otimes \mu )}{t}{(w(\cdot ,\cdot ))(t)}$$
is continuous. By [C] Proposition 6.1.4.9 a), the map
$$\mad{L^2(\mu \otimes \mu )}{\lbb{2}{L^2(\mu )}}{k}{\overbrace{k} }$$
is an isometry of Banach spaces. Since for all $t\in T$ 
$$\varphi _t\tilde{w}=\overbrace{(w(\cdot ,\cdot )(t))}$$
we get
$$(\theta _n(\tilde{w} ))(t)=\theta _n(\varphi _t\tilde{w})=\theta _n(\overbrace{w(\cdot ,\cdot )(t)})$$
for all $n\in \bn$ and so
$$\si{n\in \bn}(\theta _n(\tilde{w} ))(t)^2=\si{n\in \bn}\theta _n(\overbrace{w(\cdot ,\cdot )(t)})^2=\n{(w(\cdot ,\cdot ))(t)}_2^2.$$
Thus the map
$$\mad{T}{\br}{t}{\si{n\in \bn}}(\theta _n(\tilde{w} ))(t)^2$$
is continuous and $\tilde{w}\in \lca{2} $.\qed
 
\begin{center}
{\bfseries REFERENCES}
\end{center}
\begin{flushleft}
[C] Constantinescu, Corneliu, C*-algebras. Elsevir, 2001. \newline
[L] Lance, E. Christopher, Hilbert C*-modules, Cambridge University Press, 1995. \newline
[W] Wegge-Olsen, N. E., K-Theory and C*-Algebras, Oxford University  Press, 1993.
\end{flushleft}

\begin{center}
SUBJECT INDEX
\end{center}
\hspace{6mm} $u$-orthonormal sequence in $K$ (\dd{23.1a})

Schatten decomposition (\h{28.1a})

H\"older inequality (\cor{1.1a})

trace (\dd{1.2})

Hilbert-Schmidt operators on $H$ (\dd{1.1b})

E-$\mu $-integrable (\dd{39})

\begin{center}
SYMBOL INDEX
\end{center}
\hspace{6mm} $\psi $, $\psi _t$, $\varphi _t$ (\dd{12})

$\theta _n(u)$ (\pr{23.1})

$\xi (t)$, $U_n(u)$, $e_n(u)$ (\dd{23.1a})

$\lc{E}{p}{H}$, $\n{\cdot }_p$ (\dd{25})

$\Omega $ (\dd{31.1})

tr (\dd{1.2}) 

$ux$ (\pr{4.2})

$\int g\,\mathrm{d}\mu =\int g(s)\,\mathrm{d}\mu (s)$ (\dd{39})

\begin{flushright}
{\scriptsize \hspace{-5mm} Corneliu Constantinescu\\
Bodenacherstr. 53\\
CH 8121 Benglen\\
e-mail: constant@math.ethz.ch }
\end{flushright}
\end{document}